\documentclass{amsart}
\usepackage{latexsym}
\usepackage{amsfonts}
\usepackage{graphicx}
\usepackage{amsmath}

\newtheorem{theorem}{Theorem}[section]

\newtheorem{corollary}{Corollary}[section]

\newtheorem{definition}{Definition}[section]
\newtheorem{example}{Example}[section]

\newtheorem{lemma}{Lemma}[section]
\newtheorem{notation}{Notation}[section]

\newtheorem{proposition}{Proposition}[section]
\newtheorem{remark}{Remark}[section]

\begin{document}
\title[Uniform Closure]{Uniform Closure of Dual Banach Algebras}
\author[M. Amini]{Massoud Amini}
\address{Department of Mathematics, University of Shahid Beheshti\\
Evin, Tehran 19839, Iran
\linebreak
\indent Department of Mathematics and Statistics\\ University of Saskatchewan
\\106 Wiggins Road, Saskatoon\\ Saskatchewan, Canada S7N 5E6\\mamini@math.usask.ca}
\thanks{* Supported by a grant from Shahid Beheshti University.}
\maketitle

\begin{abstract}
We give a characterization of the ''uniform closure'' of the dual of a $%
C^{\ast }$-algebra. Some applications in harmonic analysis are given.
\end{abstract}

\maketitle

\section{Four topologies on the unit ball of a $C^{\ast }$-algebra}

Let $A$ be a $C^{\ast }$-algebra and $P(A)$,$S(A)$,$A_{1}^{\ast }$,
and $A^{\ast }$denote the\ pure state space, state space, closed
dual unit ball, and the dual of\ $A$, respectively, all equipt with
the w$^{\ast }$-topology. Let $A_{1}$denote the\ closed unit ball of 
$A$. Following [L], we consider the following four\ topologies on $%
A_{1}$:

\begin{enumerate}
\item  (w) $a_{i}\rightarrow 0$ iff $<a_{i},f>\rightarrow 0$ \ \ $(f\in
A^{\ast })$, \-

\item  (wo) $a_{i}\rightarrow 0$ iff $<a_{i}x,f>\rightarrow 0$ \ \ $(x\in
A,f\in A^{\ast })$, 

\item  (so) $a_{i}\rightarrow 0$ iff $\ ||a_{i}x||\rightarrow 0$ \ \ $(x\in
A)$, 

\item  (uc) $a_{i}\rightarrow 0$ iff $\ a_{i}\rightarrow 0$ , uniformely on w%
$^{\ast }$-compact subsets\ of $S(A)$. \ \ \ \ \ \ \ \ \ \ \
\ \ \ \ \ \ \ \ \ \ \ \ \ \ \ \ \ \ \ \ \ \ \ \ \ \ \ \ \ \ \ \ \ \ \ \ \ \
\ \ 
\end{enumerate}

Note that (so) is just the restriction of the strict topology of
the\ multiplier algebra $M(A)$ to $A_{1}$. Also note that if in (4)
one requires\ the uniform convergence on w$^{\ast }$-compact subsets
of $A_{1}^{\ast }$ (instead of\ $S(A)$ ), one get nothing but the
norm topology (Banach-Alaoglu).

\begin{lemma}
(Akemann-Glimm) Let $H$ be a Hilbert space and\ $S,T\in B(H)$. Let $%
g:\mathbf{\Bbb R}\rightarrow \mathbb{R}$ be a non negative Borel
function\ and $0<\theta <1$. Assume moreover that

\begin{enumerate}
\item  $0\leq T\leq 1$, $S=S^{\ast }$, and $S\geq T$ ,

\item  $g\geq 1$ on $[\theta ,+\infty )$ ,

\item  $<T\zeta ,\zeta >\geq 1-\theta $, for some $\zeta \in H$.

Then $<g(S)\zeta ,\zeta >\geq 1-4\sqrt{\theta }.  $
\end{enumerate}

\begin{proof}
Lemma 11.4.4 of [D].
\end{proof}
\end{lemma}

\begin{proposition}
Topologies w, wo, and so coincide on $A_{1}$ and\ they are stronger
than uc.
\end{proposition}

\begin{proof}
(wo$\subset $w). Given $f\in A^{\ast }$ and $x\in A$ , consider the\
Arens product $x.f\in A^{\ast }$ defined by $x.f(a)=f(ax)$ \ $(a\in A)$
.\ If $\{a_{i}\}\subset A_{1}$ and $a_{i}\rightarrow 0$ (w) then $%
<a_{i}x,f>=<a_{i},x.f>\rightarrow 0$ ,\ i.e. $a_{i}\rightarrow 0$
(wo).

(w$\subset $wo). By the Cohen Factorization theorem [DW] ,\ we have $%
A^{\ast }=A.A^{\ast }=\{x.f:$ \ $x\in A,f\in A^{\ast }\}$. Now if\ $%
\{a_{i}\}\subset A_{1}$ and $a_{i}\rightarrow 0$ (wo), then given $g\in
A^{\ast }$ , choose $x\in A $ $\ $and $f\in A^{\ast }$ such that $%
g=x.f$ . Then $<a_{i},g>=<a_{i}x,f>\rightarrow 0$ ,\ i.e. $%
a_{i}\rightarrow 0$ (w).

(wo$\subset $so). Trivial.

(so$\subset $wo). We adapt the proof of theorem 1 of [L]. Given\ $%
f\in A^{\ast }$ , we need only to show that if $f$ restricted to $A_{1}$
is\ so-continuous, then it is also wo-continuous. To this end,
assume\ the so-continuity and note that each $a\in A$ could be
associated\ with the (bounded) linear operator on $A$ , taking $x\in
A$ to $ax$. If\ $N=\ker (f)$, then by convexity of $A_{1}$, we
have\ $(N\cap A_{1})^{-so}\cap A_{1}=N\cap A_{1}$ (theorem 13.5 of
[KN]). Hence\ $(N\cap A_{1})^{-wo}\cap A_{1}=N\cap A_{1}$ (corollary
5 of [DS]) , and so $f$ is\ wo-continuous.

(uc$\subset $wo). We use an idea of [AAP]. Take a w$^{\ast }$%
-compact\ subset $K$ of $S(A)$ and let $f\in K$. Then, given $\theta
>0$ there exists\ $a_{f}\in A$ such that $0\leq a_{f}\leq 1$ and $%
f(a_{f})>1-\theta $ (see the proof\ of lemma 4.5 in [AAP]). Take the
w$^{\ast }$-neighbourhood\ $V_{f}=N(f,a_{f},\theta /2)=\{g\in S(A)$ $%
:|g(a_{f})-f(a_{f})|<\theta /2\}$ \ of $f$ in\ $S(A)$. Then given $%
g\in V_{f}$ we have\ $g(a_{f})\geq f(a_{f})+\theta /2>1-\theta
+\theta /2=1-\theta /2$. Cover $K$ by\ $\{V_{f}\}_{f\in K}$ and use w%
$^{\ast }$-compactness of $K\;$to get $n\geq 1$ and\ $%
f_{1},f_{2},...,f_{n}\in K$ such that $K\subset V_{f_{1}}\cup ...\cup
V_{f_{n}}$ . Put $a_{i}=a_{f_{i} }$ \ $(i=1,...,n)$ and $
a=a_{1}+...+a_{n}$ . Let $g:\Bbb R\rightarrow \mathbb{R}$ be a\
continuous function which is $0$ on $(-\infty ,0]$, $1$ on $[2\sqrt{\theta }%
,+\infty )$, and\ linear on $[0,2\sqrt{\theta }]$. Put $b=g(a)$,
then $0\leq b\leq 1$ . Now given\ $f\in K$ we have $f\in V_{f_{i}}$
, for some $i$, say $i=1$. Then\ $f(a_{1})\geq 1-\theta /2$ . On the
other hand, there is a cyclic\ representation $\{\pi ,H,\zeta \}$ of $%
A$ such that $f(x)=<\pi (x)\zeta ,\zeta > $ \ $(x\in A)$ . Take $T=\pi
(a_{1})$ and $S=\pi (a)$, then clearly\ $0\leq T\leq 1$ and $S\geq T$
. Hence by above lemma, 
\begin{eqnarray*}
f(b)=<\pi (b)\zeta ,\zeta > & = & <\pi (g(a))\zeta,\zeta > \\ 
& = & <g(\pi (a))\zeta ,\zeta > = <g(S)\zeta ,\zeta > \\
& \geq & 1-4\sqrt{\theta /2}  \geq 1-4\sqrt{\theta } .
\end{eqnarray*}
Now consider a net $\{a_{i}\}\subset A_{1}$ such that $a_{i}\rightarrow 0$
(so),\ then inside $M(A)$ we can write $a_{i}=a_{i}b+a_{i}(1-b)$%
,\ $(1-b)^{2}\leq (1-b)$, and $f(1)=||f||=1$. Hence, for each $
i $
\begin{eqnarray*}
|f(a_{i})|\leq |f(a_{i}b)|+|f(a_{i}(1-b))| 
 & \leq & ||a_{i}b||+f(a_{i}a_{i}^{\ast
})^{1/2}f((1-b)^{2})^{1/2} \\ 
& \leq & ||a_{i}b||+||a_{i}a_{i}^{\ast
}||^{1/2}f(1-b)^{1/2} \\
& \leq & ||a_{i}b||+(1-(1-4\sqrt{\theta }
))^{1/2} \\
& = & ||a_{i}b||+2\sqrt[4]{\theta } .
\end{eqnarray*}
Hence $\sup_{f\in K}|f(a_{i})|\leq ||a_{i}b||+2\sqrt[4]{\theta }$, and so $%
a_{i}\rightarrow 0 $ uniformely on $K$, as required.
\end{proof}

\begin{corollary}
If $A$ is a $C^{\ast }$-algebra, $A_{1}$ is the unit ball of $A$,
and\ $f:A\rightarrow \mathbb{C}$ is continuous with respect to (uc)
, then f is\ continuous with respect to (w). \-
\end{corollary}

\section{Dual Algebras}

\begin{notation}
(Walter) [W] If $A$ is a $C^{\ast }$-algebra, $\frak{L}(A)$, $\frak{P}(A)$%
,\ and $\frak{D}(A)$ denote the collection of all bounded,
completely\ positive, and completely bounded linear maps of $A$ into 
$A$,\ respectively. $\frak{D}(A)$ is called the \textbf{dual algebra 
}of $A$.
\end{notation}

It can be shown that $\frak{D}(A)$ is a Banach algebra with\
conjugation (this is the same as involution , except that it
preserves\ the order of multiplication) , and if $\frak{B}(A)$ is
the closed linear span\ of $\frak{P}(A)$ in $\frak{D}(A)$ (with
respect to the completely bounded norm)\ then $\frak{B}(A)\subset 
\frak{D}(A)\subset \frak{L}(A)$ [W].

\begin{definition}
(Walter) [W] Let $A$ and $B$ be Banach algebras\ with involution and
conjugation such that there are\ $C^{\ast }$-algebras $C^{\ast }(A)$
and $C^{\ast }(B)$ satisfying the following\ conditions:

\begin{enumerate}
\item  There are Banach algebra homomorphisms\ $i_{A}:A\rightarrow
C^{\ast }(A)$ and $i_{B}:B\rightarrow C^{\ast }(B)$ which are
one-one,\ onto a dense subalgebra, and preserve involution.

\item  There are norm decreasing Banach algebra\ isomomorphisms $%
j_{A}:A\rightarrow \frak{D}(C^{\ast }(B))$ and\ $j_{B}:B\rightarrow 
\frak{D}(C^{\ast }(A))$ which preserve conjugation.

Then $A$ and $B$ are called \textbf{dual algebras. }If the
involutions\ and conjugations of both algebras are isometric,
the\ duality is called \textbf{semi rigid}. If moreover both $j_{A}$
and $j_{B} $ are isometric, the duality is called\textbf{\
rigid. }
\end{enumerate}
\end{definition}

\begin{definition}
Consider the dual algebras $A$ and $B$. The duality\ is called 
\textbf{complete }if there are norm decreasing linear\ injections $%
k_{A}:C^{\ast }(A)^{\ast }\rightarrow M(C^{\ast }(B))$ and\ $%
k_{B}:C^{\ast }(B)^{\ast }\rightarrow M(C^{\ast }(A))$. Here $M$ stands for
the multiplier\ algebra.

\begin{example}
If $G$ is a locally compact group then the Fourier\ algebra $A(G)$
and the group algebra $L^{1}(G)$ are dual. Here\ we take $C^{\ast
}(A(G))=C_{0}(G)$ and $C^{\ast }(L^{1}(G))=C^{\ast }(G)$. The\
duality is rigid [W] and complete [P].
\end{example}

\begin{example}
If $A$ is $M_{n}(\Bbb C)$ with shur product and trace norm\ and $B$ is $%
M_{n}(\mathbb{C})$ with usual matrix product and $L^{1}$ norm, then\ 
$A$ and $B$ are dual and duality is rigid [W] and complete. This\ is
a special case of example 3.1 for $G$ being the principal\ transitive
groupoid on $n$ elements.
\end{example}

\begin{example}
If $A$ is the $C^{\ast }$-algebra of trace class operators on\ $\ell
^{2}$ and $B$ is the subalgebra of $M_{\infty }(\mathbb{C})$ consisting of
countably\ infinite matrices with finite $L^{1}$ norm, then $A$ and $%
B$ are dual\ and duality is semi rigid [W] and complete. This is a
special\ case of example 3.1 for $G$ being the principal
transitive\ groupoid on countably many elements.
\end{example}
\end{definition}

\begin{definition}
Consider the dual algebras $A$ and $B$. The duality\ is called 
\textbf{amenable }if there are surjective isometric\ isomorphisms $%
l_{A}:C^{\ast }(A)^{\ast }\rightarrow M(B)$ and $l_{B}:C^{\ast }(B)^{\ast
}\rightarrow M(A)$.

\begin{example}
The duality of example 2.1 is amenable iff the\ locally compact group $%
G$ is amenable [L].

\begin{proposition}
Every amenable rigid duality is complete.

\begin{proof}
Since the Banach algebra $\frak{D}(C^{\ast }(B))$ is unital, the\
isometric isomorphism $j_{A}:A\rightarrow \frak{D}(C^{\ast }(B))$ uniquely
extends\ to one on $M(A)$, still denoted by $j_{A}$ . Take $k_{B}=$ $%
l_{B}\circ j_{A}$ . Now\ $k_{A}$ is constructed similarly.
\end{proof}
\end{proposition}
\end{example}
\end{definition}

\begin{remark}
Example 2.1 showes that the convers of above\ proposition is not
true.
\end{remark}

\section{Uniform closure of dual algebras}

Consider the dual algebras $A$ and $B$. If the duality is complete\textbf{,}
then using the norm decreasing linear injection\ $k_{A}:C^{\ast
}(A)^{\ast }\rightarrow M(C^{\ast }(B))$ , one can identify $C^{\ast
}(A)^{\ast }$ with a subalgebra of $M(C^{\ast }(B))$, where of course the
norm of the\ later (which is denoted by $||.||_{u}$ ) is weaker. In
this section we want to calculate the closure of $C^{\ast }(A)^{\ast }$ in $%
M(C^{\ast }(B))$,\ which we call the \textbf{uniform closure}%
.

\begin{theorem}
Consider the dual algebras $A$ and $B$. If the duality is rigid and
complete, then the closure of $C^{\ast }(A)^{\ast } $ in $M(C^{\ast
}(B))$ consists exactly of those elements $b\in M(C^{\ast }(B))$ which
satisfy the following property:

If $\{a_{n}\}$ is any sequence in the unit ball $A_{1}$ of $A$ such that $%
<a_{n},f>\rightarrow 0$ for all $f\in C^{\ast }(A)^{\ast }$, then\ $%
<a_{n},b>\rightarrow 0$.

\begin{proof}
Assume that $b$ is in the uniform closure of $C^{\ast }(A)^{\ast }$ and $%
\{a_{n}\}$ is any sequence in the unit ball $A_{1}$ of $A $ such
that $<a_{n},f>\rightarrow 0$ for all $f\in C^{\ast }(A)^{\ast }$. Let $%
\theta >0$, and take $g\in C^{\ast }(A)^{\ast }$ such that $%
||b-g||_{u}<\theta $. Then by\ assumption, $<a_{n},g>\rightarrow 0$.
Therefore
$$
\lim \sup_{n\rightarrow \infty } |<a_{n},b>|=\lim \sup_{n\rightarrow
\infty } |<a_{n},b-g>|\leq \lim \sup_{n\rightarrow \infty }
||b-g||_{u}.||a_{n}||<\theta .
$$
Hence\ $<a_{n},b>\rightarrow 0$.

Conversely suppose that $b\in M(C^{\ast }(B))\backslash (C^{\ast }(A)^{\ast
})^{-||.||_{u}}$. Then by closed graph theorem, $b$ is not\ $w$%
-continuous on $A_{1}$, where $w=\sigma (A,C^{\ast }(A)^{\ast })$. By
corollary 1.1, $b$ is not $uc$-continuous on $A_{1},$hence there is\ $\theta
>0 $ and for each norm bounded $K\subset C^{\ast }(A)^{\ast }$ and each $%
\delta >0$, there is $a_{K,\delta }\in A_{1}$ such that
$$
|<a_{K,\delta },b>|\geq
\theta  \ , \ \ |<a_{K,\delta },f>|<\delta \, \, \, (f\in K). 
$$
Fix $w^{\ast }$-compact subset $K\subset C^{\ast }(A)^{\ast }$ and put $%
a_{1}=a_{K,1}$. Then take 
$$
K_{1}=\{f\in C^{\ast }(A)^{\ast }:  \ |<a_{1},f>|\geq 1\}
$$ 
and put $a_{2}=a_{K_{1},1}$. Continuing this way we put 
$$
K_{n}=\{f\in C^{\ast }(A)^{\ast }\ :|<a_{i},f>|\geq 1/n \, (1\leq i\leq n)\}
$$ 
and $a_{n+1}=a_{K_{n},1/n}$ . Then $
<a_{n},f>\rightarrow 0$ for all $f\in C^{\ast }(A)^{\ast }$ (for those $
f $ which belong to $\cup _{n\geq 1}K_{n}$ use the defining property
of $a_{K,\delta }$'s and for others use the defining property of 
$K_{n}$'s)\ but \ $|<a_{n},b>|\geq \theta $ \ $(n\geq 1)$, and we are
done.
\end{proof}
\end{theorem}

It is clear from the proof of the above theorem that we only need to assume
a ''one way duality'' relation betwen\ two algebras. More precisely
it is enough that $A$ and $B$ satisfy the following definition.

\begin{definition}
Let $A$ and $B$ be Banach algebras with involution such that there are $%
C^{\ast }$-algebras $C^{\ast }(A)$ and $C^{\ast }(B) $ satisfying
the following conditions:

\begin{enumerate}
\item  There are Banach algebra homomorphisms $i_{A}:A\rightarrow C^{\ast
}(A)$ and $i_{B}:B\rightarrow C^{\ast }(B)$ which are one-one, onto
a\ dense subalgebra, and preserve involution.

\item  There is norm decreasing linear injection $k_{A}:C^{\ast }(A)^{\ast
}\rightarrow M(C^{\ast }(B))$. Here $M$ stands for the multiplier\
algebra.
\end{enumerate}
\end{definition}

\textit{Then }$A$\textit{\ is called {\bf semi dual} to }$B$.

\begin{example}
If $G$ is a topological groupoid with Haar system then the Fourier algebra $%
A(G)$ [RW], [R] is\ semi dual to the convolution algebra $C_{c}(G)$
(except that here $C_{c}(G)$ is only a normed *-algebra). Here we
take\ $C^{\ast }(A(G))=C_{0}(G)$ and $C^{\ast }(C_{c}(G))=C^{\ast
}(G)$.
\end{example}

\begin{example}
If $S$ is a commutative separative foundation topological semigroup then the 
$L^{\infty }$-representation\ algebra $R(S)\;$on $S$ [L] is semi dual
to the semigroup algebra $M(S)$. Here we take $C^{\ast }(R(S))=C_{b}(S)$
and\ $C^{\ast }(M(S))=W^{\ast }(S)$ [La].
\end{example}

\begin{example}
If S is a foundation topological *-semigroup whose *-representations
separate the points of S, then the Fourier algebra $A(S)$ [AM] is semi dual
to the\ semigroup algebra $M_{a}(S)$. Here we take $C^{\ast
}(A(S))=C_{0}(S)$ and $C^{\ast }(M_{a}(S))=C^{\ast }(S)$ [AM].
\end{example}

\begin{example}
If $S$ is an inverse semigroup then the Fourier algebra $A(S)$ [AM2] is
semi dual to the semigroup\ algebra $\zeta ^{1}(S)$. Here we take $%
C^{\ast }(A(S))=C_{0}(S)$ and $C^{\ast }(\zeta ^{1}(S))=C^{\ast }(S)$
[AM2].
\end{example}

\begin{theorem}
Consider the involutive Banach (normed) algebras $A$ and $B$. If $A$ is semi
dual to $B$, then the closure\ of $C^{\ast }(A)^{\ast }$ in $%
M(C^{\ast }(B))$ consists exactly of those elements $b\in M(C^{\ast }(B))$
which satisfy the following property:

If $\{a_{n}\}$ is any sequence in the unit ball $A_{1}$ of $A$ such that $%
<a_{n},f>\rightarrow 0$ for all $f\in C^{\ast }(A)^{\ast }$, then\ $%
<a_{n},b>\rightarrow 0$.
\end{theorem}

As a corollary we get the main result of [L], whose proof is adapted in our
theorem 3.1.

\begin{corollary}
\textbf{[L] }Let $S$ be a commutative separative foundation semigroup with
identity and let $R(S)$ denote the\ $L^{\infty }$-representation
algebra on $S$. Then for a function $f\in C_{b}(S)$ the following are
equivalent:

\begin{enumerate}
\item  $f\in R(S)^{-||.||_{\infty }}$.

\item  If $\{\mu_{n}\}$ is any sequence in the unit ball of $M_{a}(S)$ such
that $\mu_{n}^{\symbol{94}}(\chi)\rightarrow 0$ as $n\rightarrow \infty $, for
all $\chi\in S^{\symbol{94}\text{, }}$then $\int_{S} $ $%
fd\mu_{n}\rightarrow 0$, as $n\rightarrow \infty $.

\begin{proof}
Example 3.2 and theorem 3.2.
\end{proof}
\end{enumerate}
\end{corollary}

As another corollary we get the non commutative version of the above result,
which has been our main motivation\ to write this paper.

\begin{corollary}
Let $S$ be a foundation topological *-semigroup with identity whose
*-representations separate the points of S, and let $B(S)$ denote the
Fourier-Stieltjes algebra on\ $S$. Then for a function $f\in
C_{b}(S) $ the following are equivalent:

\begin{enumerate}
\item  $f\in B(S)^{-||.||_{\infty }}$.

\item  If $\{\mu_{n}\}$ is any sequence in the unit ball of $M_{a}(S)$ such
that $\int_{S}$ $gd\mu_{n}\rightarrow 0$ as $n\rightarrow \infty $, for all $%
g\in P(S)^{\text{, }}$then $\int_{S} $ $fd\mu_{n}\rightarrow 0$, as $%
n\rightarrow \infty $.

\begin{proof}
Example 3.3 and theorem 3.2.
\end{proof}
\end{enumerate}
\end{corollary}

As far as I know, this result\  is new even for locally
compact groups (although $B(G)^{-||.||_{\infty }}$ has been studied in other
directions (see for instance [C]).

\begin{corollary}
Let $G$ be a topological group and $m$ be a left Haar measure on $G$, and let $B(G)$ denote the Fourier-Stieltjes
algebra on $G$. Then for a\ function $f\in C_{b}(G)$ the following
are equivalent:

\begin{enumerate}
\item  $f\in B(G)^{-||.||_{\infty }}$.

\item  If $\{f_{n}\}$ is any sequence in the unit ball of $L^{1}(G)$ such
that $\int_{G}$ $gf_{n}dm\rightarrow 0$ as $n\rightarrow \infty $, for all $%
g\in P(G)^{\text{, }}$then $\int_{G} $ $ff_{n}dm\rightarrow 0$, as $%
n\rightarrow \infty $.
\end{enumerate}
\end{corollary}

\textbf{References}

[AAP] C.A. Akemann, J. Anderson, G.K. Pedersen, \textbf{Diffuse sequences and} 

\indent \, \, \, 
\textbf{perfect $C^*$-algebras}, Trans. Amer. Math. Soc. 298 (2)(1986), 747-762.

[AM]  M. Amini, A.R. Medghalchi, \textbf{Fourier algebra on
topological *-semi-}

\indent \, \, \, 
\textbf{groups}, preprint, 2000.

[AM2]\ M. Amini, A.R. Medghalchi, Fourier algebra on inverse semigroups, 

\indent \, \, \, 
preprint, 2000.

[C] C.\ Chou, \textbf{Weakly almost periodic functions and Fourier-Stieltjes}

\indent \, \, \, 
\textbf{algebra of locally compact groups} , Trans. Amer. Math. Soc. 274

\indent \, \, \,
(1982) 141-157

[D] J.  Dixmier, C*-algebras, North Holland, 1977.

[DS]  N. Dunford, T. Schwartz , Linear operators I, Interscience,
New York, 1958.

[DW]\ R.S. Doran, J. Wichmann, Approximate identities and
factorization in 

\indent \, \, \, 
Banach modules , Lecture Notes in Mathemat ics 768,
Springer Verlag, 

\indent \, \, \, 
Berlin, 1979.

[KN]  J.L.\ Kelley, I. Namioka, Linear topological spaces, Van
Nostrand, 

\indent \, \, \, 
Princeton, 1963.

[L] M . Lashkarizadeh Bami, \textbf{On some sup-norm closure of the} 
$L^{\infty }$\textbf{-}

\indent \, \, \, 
\textbf{representation algebra R(S) of a foundation
semigroup\ S}, semigroup 

\indent \, \, \, 
forum 52 (1996) 389-392.

[P] G.K. Pedersen, C*-algebras and their authomorphism groups, Academic

\indent \, \, \, 
Press, 1979.

[R] J. Renault, \textbf{The Fourier algebra of a measured groupoid and its}

\indent \, \, \, 
\textbf{multipliers}, J. Funct. An al. 145 (1997) 455-490.

[RW] A.  Ramsay, M.E. Walter, \textbf{Fourier-Stieltjes algebra of}
locally compact 

\indent \, \, \, 
\textbf{groupoids}, J. Funct. Anal. 148 (1997) 314-367.

[W]  M.E. Walter, \textbf{Dual algebras}, Math.\ Scand. 58 (1986)
77-104.

\end{document}